\documentclass[reqno]{amsart}

\usepackage{amsfonts}
\usepackage{amscd}
\usepackage{amsbsy}
\usepackage{amsxtra}
\usepackage{amssymb}
\usepackage{epsfig}
\usepackage{epic,eepic}
\usepackage{graphicx}
\usepackage[all]{xy}
\usepackage{xypic}
\usepackage{psfrag}
\usepackage{amsmath}
\usepackage{amsthm}
\usepackage{setspace}
\usepackage{hyperref}
\usepackage{url}
\usepackage{here}
\usepackage{todonotes}
\usepackage{dsfont}
\newlength{\halfbls}
\usepackage{tikz}
\usetikzlibrary{calc}
\usetikzlibrary{decorations.pathreplacing,decorations.markings}
\usetikzlibrary{arrows}

\DeclareRobustCommand{\SkipTocEntry}[9]{}

   \newcommand{\CC}{\mathbb{C}}
  \newcommand{\HH}{\mathbb{H}}

\newcommand{\QQ}{\mathbb{Q}} \newcommand{\RR}{\mathbb{R}}  

\newcommand{\ZZ}{\mathbb{Z}}

    \newcommand{\Aut}{{\rm Aut}}
\newcommand{\Pic}{{\rm Pic}}  

\newcommand{\ord}{{\rm ord}}

\newcommand{\rank}{{\rm rank}}

      \def\frakb{\mathfrak b}

\newcommand{\frakd}{\mathfrak{d}}

\newcommand{\cLL}{{\mathcal L}} 
 \newcommand{\cFF}{{\mathcal F}}  
\newcommand{\cDD}{{\mathcal D}}   
\newcommand{\cCC}{{\mathcal C}}

 \newcommand{\cOO}{{\mathcal O}} 
\newcommand{\calO}{{\mathcal O}}
\newcommand{\calF}{{\mathcal F}}
\newcommand{\cPP}{{\mathcal P}}

\newcommand{\SO}{{\rm SO}}
\newcommand{\Orth}{\operatorname{O}}

\newcommand{\norm}{\operatorname{N}}
\newcommand{\Mp}{{\rm Mp}}

\newcommand{\wOL}{\widetilde{\Orth}^+\!(L)}

    \newcommand{\ve}{{\varepsilon}}

\newtheorem{Defi}{Definition}[section]  
    \newtheorem{Prop}[Defi]{Proposition}
\newtheorem{Lemma}[Defi]{Lemma}    \newtheorem{Cor}[Defi]{Corollary}
\newtheorem{Thm}[Defi]{Theorem} 
\newtheorem*{Thm*}{Theorem}

\newcommand{\leg}[2]{\left( \frac{#1}{#2} \right)}
\newcommand{\zxz}[4]{\begin{pmatrix} #1 & #2 \\ #3 & #4 \end{pmatrix}}

\newcommand{\kzxz}[4]{\left(\begin{smallmatrix} #1 & #2 \\ #3 & #4\end{smallmatrix}\right) }

\def\={\;=\;}  \def\+{\,+\,}

\newcommand\R{\mathbb R}  \newcommand\C{\mathbb C}  \newcommand\Z{\mathbb Z}  \newcommand\Q{\mathbb Q} 
     
\renewcommand\H{\mathbb H}      

       \def\p{\varphi}   \def\2{\pi_2}

\def\SL{{\rm SL}}

\def\mat#1#2#3#4{\begin{pmatrix}#1&#2\\#3&#4\\ \end{pmatrix}}       

\def\be{\begin{equation}}   \def\ee{\end{equation}}     \def\bes{\begin{equation*}}    \def\ees{\end{equation*}}
\def\ba{\be\begin{aligned}} \def\ea{\end{aligned}\ee}   \def\bas{\bes\begin{aligned}}  \def\eas{\end{aligned}\ees}

\newcounter{savedtocdepth}
\newcommand*{\SaveTocDepth}[1]{%
  \addtocontents{toc}{%
    \protect\setcounter{savedtocdepth}{\protect\value{tocdepth}}%
    \protect\setcounter{tocdepth}{#1}%
  }%
}

\title[Cones of Heegner divisors]
{Cones of Heegner divisors}

\author[Jan H.~Bruinier]{Jan
Hendrik Bruinier}
\address{Fachbereich Mathematik,
Technische Universit\"at Darmstadt, Schlossgartenstrasse 7, D--64289
Darmstadt, Germany}
\email{bruinier@mathematik.tu-darmstadt.de}

\author{Martin M\"oller}
\address{
Institut f\"ur Mathematik, Goethe--Universit\"at Frankfurt,
Robert-Mayer-Str. 6--8,
60325 Frankfurt am Main, Germany\\
}
\email{moeller@math.uni-frankfurt.de}


\thanks{The first author is partially supported by DFG grant BR-2163/4-2.}

\begin{document}
\bibliographystyle{halpha}

\begin{abstract}
We show that the cone of primitive Heegner divisors
is finitely generated for many orthogonal Shimura varieties, 
including the moduli space of polarized K3-surfaces.
The proof relies on the growth of coefficients of
modular forms.
\end{abstract}
\maketitle

\tableofcontents

\noindent
\SaveTocDepth{1}

\section{Introduction}

Let $X$ be a projective algebraic variety. Both the pseudo-effective
cone of~$X$ (the closure of the cone ${\rm Eff}(X)$ of effective divisors)
and  dually (by \cite{BDPP}) the cone of movable curves are important
geometric invariants of~$X$ that are notoriously hard to compute.
The same claim can be made for the cone of curves, or dually for
the cone of nef divisors. Mori's cone theorem is the most important
result implying that extremal rays in the cone of curves in
the half-space where $K_X$ is negative can accumulate at most
towards the plane $K_X^\perp$. On the other hand, 
abelian surfaces provide examples where the all these cones are round.
For varieties of general type the knowledge about these cones is very 
limited and the existing polyhedrality results (for example \cite{KKL}) 
have hypothesis that are restrictive or hard to verify.
E.g.\ it was shown recently by Mullane (\cite{Mullane}) that
the effective cone ${\rm Eff}(M_{g,n})$ of the moduli space of curves
is not finitely generated for $g \geq 2$ and $n \geq g+1$.
\par
In this note we deduce from properties of modular froms the
polyhedrality of a natural subcone of the pseudo-effective cone
on orthogonal Shimura varieties. For concreteness, we work in the introduction
with the moduli space $\cFF_{2d}$ of polarized $K3$-surfaces of degree~$2d$.
For every integer~$d$ this space is the quotient 
$\cFF_{2d} = \wOL \backslash \cDD_{2d}$ of a $19$-dimensional complex 
domain $\cDD_{2d}$ by an arithmetic 
lattice~$\wOL$, see Section~\ref{sec:background} for the background on 
lattices and this quotient.
These moduli spaces carry an infinite collection of divisors, the 
Noether-Lefschetz divisors, geometrically defined
by as the loci of $K3$-surfaces with Picard group of rank~$\geq 2$.
These divisors are also called Heegner divisors.
They are not irreducible in general. We recall the structure of the
irreducible components, called primitive Heegner divisors, in 
Section~\ref{sec:IrrNL}.
\par
The structure of the Picard group $\Pic(\cFF_{2d})$ 
up to torsion 
is now completely understood. It is shown in \cite{BLMM} to be
generated by Noether-Lefschetz divisors and the rank has been computed
in \cite{Bruinier_Rank}. A next step towards understanding the geometry
of $\cFF_{2d}$ would be to compute the natural cones in $\Pic(\cFF_{2d})$, 
the ample cone and the pseudo-effective cone. 
The latter contains as a subcone the cone  ${\rm Eff}^{{\rm NL}}(\cFF_{2d})$
generated by the  primitive Heegner divisors. The theorem of \cite{BLMM} 
implies in this language that ${\rm Eff}^{{\rm NL}}(\cFF_{2d})$ is 
full-dimensional, but it does not imply that the two cones coincide.
The following theorem answers a question raised by Peterson in 
\cite{PetersonK3}.
\par
\begin{Thm} \label{thm:introEff}
The cone  ${\rm Eff}^{{\rm NL}}(\cFF_{2d})$ generated by the primitive Heegner 
divisors is rational polyhedral. In particular it is finitely generated.
\end{Thm}
\par
The methods can be applied to the more general situation where the
lattice~$L$ is an even lattice of signature $(b^+,2)$ with $b^+ \geq 3$.
If~$L$ splits off a hyperbolic plane we show in Theorem~\ref{thm:redHeeg}
that the cone of Heegner divisors on the hermitian symmetric space~$\cFF_{L}$ 
associated with~$L$ is rational polyhedral. 
\par
If the lattice~$L$ splits off two hyperbolic planes, we show in
Theorem~\ref{thm:EffHeeg} generalizing Theorem~\ref{thm:introEff}
that the cone ${\rm Eff}^{{\rm H}}(\cFF_{L})$ generated by primitive Heegner
divisors is rational polyhedral. 
\par
To some extent the results carry over to the case $b^+ = 2$ of 
Hilbert modular surfaces. The details are summarized in Section~\ref{sec:HMS}.
\par
\medskip
\paragraph{\bf Acknowledgments}

We thank Alex K\"uronya for several useful discussions and comments.

\section{Vector-valued modular forms 
and orthogonal Shimura varieties}
\label{sec:background}

This section gathers background material and notation on modular forms 
associated with a lattice and orthogonal Shimura varieties.
\par
\medskip
\paragraph{\bf Vector-valued modular forms associated to a lattice}

Suppose that $(L,(\cdot,\cdot))$ is an even lattice of signature $(b^+,b^-)$ 
with quadratic form $Q(x) = \tfrac12 (x,x)$. 
We denote by $L^\vee$
the dual lattice of~$L$ and by  $D_L = L^\vee/L$ the discriminant lattice.
We let $N$ be the level of~$L$, i.e.\ the smallest integer such that
$NQ(\mu) \in \ZZ$ for every $\mu \in L^\vee$. 
We denote by $(\chi_\mu)_{\mu \in D_L}$ the standard basis
of $\CC[D_L]$.
Recall that there 
is a Weil representation $\rho_L$ of the metaplectic group $\Mp_2(\ZZ)$ on 
the group ring $\CC[D_L]$ (see e.g.\ \cite{BorcherdsSing}, 
\cite{Br-habil}). 
\par
A {\em vector-valued modular form of weight~$k \in \tfrac12 \ZZ$} 
for $\Mp_2(\ZZ)$ and the Weil representation $\rho_L$ is a holomorphic
function $f: \HH \to \CC[D_L]$ that is holomorphic at $\infty$ and that
satisfies the transformation law
$$ f(\gamma \tau) \= \sigma(\tau)^k \rho_L(\gamma) f(\tau)$$
for all $\gamma = (g,\sigma) \in \Mp_2(\ZZ)$.
This space of modular forms is denoted by $M_{k,L}$.
Such a modular form has a Fourier expansion
$$ f(\tau) \= \sum_{\mu \in D_L} \sum_{0 \leq m \in Q(\mu) + \ZZ} a_{m,\mu} \,q^m \,\chi_\mu, 
\quad \text{where} \quad  q = e^{2\pi i\tau}\,. $$
\par
We suppose throughout that $2k \equiv b^+-b^- \pmod{4}$.
Let $\Gamma_\infty' \subset \Mp_2(\ZZ)$ be the stabilizer of the cusp $\infty$.
For every half-integer $k>2$ the Eisenstein series 
\be \label{eq:def:ES}
E_{k,L}(\tau) \= \sum_{(g,\sigma) \in \Gamma_\infty' \backslash \Mp_2(\ZZ)}
\sigma(\tau)^{-2k} \cdot \rho_L(g,\sigma)^{-1} \chi_0
\ee
is a vector-valued modular form in $M_{k,L}$.
There are two subspaces 
$$ S_{k,L}\subset M_{k,L}^0
\subset M_{k,L}$$ 
of the space of modular forms that are relevant in the sequel, 
namely the subspace of cusp forms $ S_{k,L}$ and the intermediate
space $M_{k,L}^0$ of forms whose constant term is supported at
the trivial element of $D_L$ only. Following \cite{PetersonK3} we refer to this
space as {\em almost cusp forms}. Obviously
\be \label{eq:ACdecomp}
M_{k,L}^0 \=  S_{k,L} \oplus \langle E_{k,L} \rangle\,.
\ee
\par
\medskip
\paragraph{\bf Eisenstein series}
The coefficients of the Fourier expansion of the Eisenstein series 
\be \label{eq:coeffES}
E_{k,L}(\tau) \= \sum_{\mu \in D_L} \sum_{m \geq 0} e_{k,L}(m,\mu) q^m \chi_\mu
\ee
have been computed explicitly in~\cite{BrKuss}. 
First, the $e_{k,L}$ are rational numbers and 
(compare \cite[Proposition~2.1]{BrPrescribed} or \cite[Theorem~4.6]{BrKuss})
\be \label{eq:signEIS}
(-1)^{(2k-b^+ + b^-)/4}\, e_{k,L}(m,\mu) \, \geq \, 0\,.
\ee
The constant term of $E_{k,L}$ is given by $\chi_0\in \C[D_L]$.
\par
{\em We now specialize to the case $k = (b^+ + b^-)/2 = {\rm rank}(L)/2$,}
which is the relevant case for the geometric application.  In this case, the 
coefficients of Fourier expansion of the Eisenstein series are given by
the following formulas (see \cite[Theorem~4.6]{BrKuss}). For a discriminant
$D \in \ZZ \setminus\{0\}$ we define the Dirichlet character 
$\chi_D = \leg{D}{a}$ and the divisor sums with character
$$\sigma_s(a,\chi) \= \sum_{d\mid a} \chi(d) d^s \,. $$
For $\mu\in D_L$ we let 
$d_\mu=\min\{b\in \Z_{>0}:\; b\mu=0\}$ be the order of $\mu$. For 
$m\in \Z+Q(\mu)$ we denote by $N^L_{m,\mu}(a)$ the mod-$a$ representation number
\[
N^L_{m,\mu}(a)= \left|\{ r\in L/aL:\; Q(r+\mu)\equiv m\pmod{a}\}\right|.
\]
We will frequently drop the superscript $L$ if it is clear from the context. Moreover, we
define 
\bes
w_p\=w_p(m,\mu)\= 1+2\ord_p(2d_\mu m)\,.
\ees
The reason for introducing this is that the normalized representation nunbers
$p^{\nu(1-2k)} N_{m,\mu}(p^{\nu})$ are independent of~$\nu$ if $\nu \geq w_p$
(see Hilfssatz~13 in \cite{siegel35}).
\par
{\em Suppose that $b^++b^-$ is even.} Then, for $\mu \in D_L$ 
and $m\in \Z+Q(\mu)$ positive, the coefficients are  given by
\ba \label{eq:EISeven}
e_{k,L}(m,\mu) 
&= \frac{(2\pi)^{k} m^{k-1}(-1)^{b^-/2}}{\sqrt{|L'/L|}\Gamma(k)}
\cdot \frac{1}{L(k,\chi_{4D})}\cdot \ve_{m,\mu}
\ea
where $D$ denotes the discriminant $D=(-1)^{(b^++b^-)/2}|D_L|$ and where
\begin{flalign} 
\ve_{m,\mu} &\= \prod_{\substack{\text{$p$ prime}\\p\mid 2m N}}  
\frac{p^{w_p(1-2k)} N_{m,\mu}(p^{w_p})}
{1 -\chi_{D}(p) p^{-k}} \label{eq:epseven1} \\
&\= \sigma_{1-k}(d_\mu^2m,\chi_{4D}) 
\prod_{\substack{\text{$p$ prime}\\p\mid 2N}} \frac{N_{m,\mu}(p^{w_p})}{p^{(2k-1)w_p}}.
\label{eq:epseven2}
\end{flalign}
\par
\smallskip
{\em Suppose that  $b^++b^-$ is odd.} We write $m d_\mu^2= m_0 f^2$ for 
positive integers~$m_0$ and~$f$ with $(f,2N)=1$ and $\ord_p(m_0)\in \{0,1 \}$ for all primes 
$p$ coprime to $2N$. In this case 
\be \label{eq:EISodd} 
e_{k,L}(m,\mu) \= \frac{(2\pi)^{k} m^{k-1}(-1)^{b^-/2}}{\sqrt{|L'/L|}\Gamma(k)}
\cdot \frac{L(k-1/2,\chi_{D'})}{\zeta(2k-1)} \cdot \ve_{m,\mu}
\ee
where $D'$ now denotes the discriminant $D'=2(-1)^{(b^++b^-+1)/2}m_0|D_L|$ 
and $\mu(\cdot)$ the M\"obius function, and where 
\begin{flalign} 
\ve_{m,\mu} &\= \prod_{\substack{\text{$p$ prime}\\p\mid 2mN}} \frac{ 1 -\chi_{D'}(p) 
p^{1/2-k}}{1- p^{1-2k}}
p^{w_p(1-2k)} N_{m,\mu}(p^{w_p})   \label{eq:epsod1}\\ 
&\=  \sum_{g\mid f} \mu(g)\chi_{D'}(g)g^{1/2-k}\sigma_{2-2k}(f/g)
\prod_{\substack{\text{$p$ prime}\\p\mid 2N}} \frac{N_{m,\mu}(p^{w_p})}
{(1-p^{1-2k})p^{(2k-1)w_p}}\,.  \label{eq:epsod2}
\end{flalign}
\par
\medskip
\paragraph{\bf Estimates for representation numbers}
The following lemmas are used to give lower and upper bounds for the coefficients~$e_{k,L}$.
Note that the representation numbers $N_{m,\mu}(a)$ are weakly multiplicative
in~$a$. 
\par

We let $U$ be the even unimodular lattice of signature $(1,1)$, realized as $\Z^2$ with the quadratic form $Q((x,y))=xy$ (a hyperbolic plane). The following lemma is well known.

\begin{Lemma} \label{lem:NU}
Let $m\in \Z$ and $\nu\in \Z_{\geq 0}$. Then
\[
N^U_{m,0}(p^\nu) = \begin{cases} 
(\ord_p(m) +1) (1-1/p)p^\nu ,&\text{if $\ord_p(m)<\nu$,}\\
\nu (1-1/p)p^\nu+p^\nu ,&\text{if $\ord_p(m)\geq \nu$.}
\end{cases}
\]  
\end{Lemma}

\begin{Cor}
\label{cor:estU}
We have
\[
1-1/p \leq p^{-\nu} N^U_{m,0}(p^\nu) \leq \nu +1 .
\]
\end{Cor}

\begin{Lemma} \label{le:RepLowBd}
Suppose that $L = L_1 \oplus U$ for an even lattice $L_1$ 
with 
$\rank(L_1)=\rank(L)-2$. Let $\mu\in D_L$ and $m\in \Z+Q(\mu)$.
Then 
\be \label{eq:RepLow}
1-1/p\leq  p^{(1-2k)\nu}N^L_{m,\mu}(p^{\nu} ) \leq   \nu+1\,.
\ee  
\end{Lemma}
\par
\begin{proof} 
Any lattice element~$\lambda \in L$ can be written  
in the basis
of $L_1 \oplus U$ as $\lambda = (\lambda_1,x,y)$ with 
$\lambda_1\in L_1$ and $x,y\in \Z$.
Then $Q(\lambda) = Q(\lambda_1) + xy$, and we may suppose that $\mu = (\mu_1,0,0)$ 
since $U$ is self-dual. 
We have
\begin{align*}
N^L_{m,\mu}(p^{\nu} ) &= |\{ \lambda_1\in L_1/p^\nu L_1, \,(x,y)\in U/p^\nu U\mid \; Q(\lambda_1+\mu_1) +xy\equiv m \pmod{p^\nu}\}|\\
&= \sum_{\lambda_1\in L_1/p^\nu L_1} N^U_{m-Q(\lambda_1+\mu_1),0}(p^\nu).
\end{align*}
Hence the claimed bounds follow directly from Corollary 
\ref{cor:estU}.
\end{proof}
\par
We now derive similar results for lattices which split two hyperbolic planes over $\Z$. 
It turns out that we get slightly stronger bounds in this case.
\par
\begin{Lemma} \label{lem:NU2}
Let $m\in \Z$ and $\nu\in \Z_{\geq 0}$. Then
\[
N^{U\oplus U}_{m,0}(p^\nu) = \begin{cases} 
p^{3\nu}(1+p^{-1})(1-p^{-\ord_p(m)-1}) ,&\text{if $\ord_p(m)<\nu$,}\\
p^{3\nu} (1+p^{-1}-p^{-\nu-1}) ,&\text{if $\ord_p(m)\geq \nu$.}
\end{cases}
\]  
\end{Lemma}
\par
\begin{proof}
For the proof we briefly  put  $M=U\oplus U$.
The  statement for $\ord_p(m)<\nu$ follows from a result of Siegel, see e.g.~\cite[Theorem~4.5]{BrKuss}. The second statement can be deduced from the first, since
\begin{align*}
N^{M}_{0,0}(p^\nu)&=\#(M/p^\nu M)-\sum_{\substack{a\in \Z/p^\nu\Z\\ a\neq 0}} 
N^{M}_{a,0}(p^\nu)\\
&=p^{4\nu}-\sum_{j=0}^{\nu-1} \sum_{b\in (\Z/p^{\nu-j}\Z)^\times}  
N^{M}_{p^j b,0}(p^\nu)\\
&=p^{4\nu}-\sum_{j=0}^{\nu-1} (p^{\nu-j}-p^{\nu-j-1}) p^{3\nu}(1+p^{-1})(1-p^{-j-1}). 
\end{align*}
Computing the latter sum, we obtain the assertion.
\end{proof}

\begin{Cor}
\label{cor:estU2}
We have
\[
1-p^{-2} \leq p^{-3\nu} N^{U\oplus U}_{m,0}(p^\nu) \leq  1+p^{-1} .
\]
\end{Cor}

\begin{Lemma} \label{le:RepLowBd2}
Let $L$ be an even lattice of rank $2k=b^++b^-$. Suppose that $L = L_0 \oplus 
U\oplus U$ for an even lattice $L_0$ of rank $2k-4$.
Let $\mu\in D_L$ and $m\in \Z+Q(\mu)$.
Then for all primes $p$ and all $\nu\in \Z_{\geq 0}$ we have 
\be \label{eq:RepLow2}
1-p^{-2}\leq  p^{(1-2k)\nu}N^L_{m,\mu}(p^{\nu} ) \leq   1+p^{-1} \,.
\ee  
\end{Lemma}
\par
\begin{proof} 
Any lattice element~$\lambda \in L$ can be uniquely written as 
$\lambda=\lambda_0+\lambda_1$ with $\lambda_0\in L_0$ and $\lambda_1\in U
\oplus U$. Then  $Q(\lambda)=Q(\lambda_0)+Q(\lambda_1)$ and 
\begin{align*}
N^L_{m,\mu}(p^{\nu} ) 
&= \sum_{\lambda_0\in L_0/p^\nu L_0} N^{U\oplus U}_{m-Q(\lambda_0+\mu),0}(p^\nu).
\end{align*}
Hence the claimed bounds follow directly from Corollary 
\ref{cor:estU}.
%
\end{proof}
\par
\medskip
\paragraph{\bf The period domain  of orthogonal Shimura varieties and Heegner divisors} 
Let $L$ be an even lattice of signature $(b^+,2)$. The Hermitian symmetric
domain $\cDD_{L}$ of the orthogonal group of this lattice can be realized as  one of the two connected components of 
$$ \cDD_{L} \cup \overline{\cDD_{L}} \= \{ z\in L_\C\,:\, Q(z) = 0 
\, \,\text{and}\,\, (z,\overline{z}) < 0\}/\C^\times\,.$$
Following \cite{GHS-2007-Inventiones}, we let $\Orth^+(L)$ be the index two subgroup 
of the orthogonal group~$\Orth(L)$ which preserves the components, that is, the subgroup of elements 
of~$\Orth(L)$ of positive spinor norm.
We let $\widetilde{\Orth}(L)$ be the discriminant kernel of $\Orth(L)$, that is, 
the kernel of the natural homomorphism $\Orth(L)\to \Aut(D_L)$, and we put
\[
\wOL  = \widetilde{\Orth}(L) \cap \Orth^+(L).
\] 
The moduli spaces we 
are interested in are the locally symmetric spaces 
 $$ \cFF_{L}(\Gamma)= \Gamma\backslash \cDD_L \quad \text{for} \quad
\Gamma \subseteq \wOL$$
a subgroup of finite index. We abbreviate $\cFF_L= \cFF_L(\wOL)$.
\par
For any vector $v \in L^\vee$ with $Q(v)>0$ the Heegner divisor $H_v \subset \cDD_{L}$
consists of the points~$z$ orthogonal to~$v$. 
For $\mu\in D_L$ and $m\in Q(\mu)+\Z$ positive, 
the group~$\wOL$ acts on 
vectors in~$\mu +L$ of norm $m$ with finitely many orbits.
Consequently, for any $\Gamma \subseteq \wOL$ the 
{\em (reducible) Heegner divisors}, defined as
$$
H_{m,\mu} \= \Gamma \,\,\backslash \,\,\Biggl( 
\sum_{\substack{v \in \mu+L
\\ Q(v) = m}}  H_v  \Biggr) \,,
$$
are well-defined in~$\cFF_{L}(\Gamma)$. These are in general neither reduced
nor irreducible. In particular for $\Gamma \subseteq \wOL$ of large index, 
$H_{m,\mu}$ may have many components. Moreover, all the components have multiplicity 
two if $\mu = -\mu$ and they all have multiplicity one otherwise 
(see Lemma~\ref{le:IrrHeeg}
below). We will discuss the passage to irreducible components in 
Section~\ref{sec:IrrNL}.
\par
The tautological line bundle $\cOO(-1)$ on $\cDD_{L}$ descends
to a line bundle~$\lambda$ on $\cFF_L$, called the {\em Hodge bundle}.
(Thus $\lambda$ is anti-ample in our notation). The Hodge bundle
plays no role in our calculation, but the intersection numbers with~$\lambda$
arise as coefficient extraction functionals, similar to the
intersection numbers with Heegner divisors, see the proof 
of Theorem~\ref{thm:redHeeg}.
\par
\medskip
\paragraph{\bf Moduli spaces of K3 surfaces and 
Noether-Lefschetz divisors } 
In the special case of the lattice $L = L_{2d}$ of signature $(19,2)$ given by  
\be \label{eq:K3L}
L_{2d} \= \langle 2d \rangle \oplus U^{\oplus 2} \oplus E_8^{\oplus 2}
\ee  
the discriminant group is $D_L \cong \ZZ/2d$. The modular variety 
$$\cFF_{2d} \cong \wOL 
\backslash \cDD_{L}$$
is closely related to  the moduli space of $2d$-polarized $K3$-surfaces. 
More precisely, 
an open subset $\cFF_{2d}^\circ$ of~$\cFF_{2d}$, the complement of
some Heegner divisors, is the moduli space of polarized K3-surfaces.
See \cite{PSSTorelli} and \cite{Mo} for the semi-polarized case
and see e.g.\ \cite{GHSSurvey} for a survey.
\par
Note that in our definition~$L$ is isomorphic to the orthogonal
complement of the polarization class~$H$ with $Q(H) = d$ in middle cohomology 
lattice
$$ L_{K3} \=  U^{\oplus 3} \oplus E_8^{\oplus 2} $$
of the K3-surface with the {\em negative} of the intersection pairing.
We let $\omega$ be a fixed generator of the first
summand of~$L$ in \eqref{eq:K3L}. Hence,  $D_L$ is generated by $\omega/2d$. 
\par
The generic algebraic $K3$-surface has a Picard group of rank one.
The (reducible) {\em Noether-Lefschetz divisors $D_{h,a}$} are the closures 
in $\cFF_{2d}$ of
the loci where the Picard group of the polarized $K3$-surfaces $(S,H)$ have
a class~$\beta$ not in the linear space of~$H$ with $Q(\beta) = h-1$
and $(\beta,H) = a$. We may assume that $0 \leq a < 2d$. 
These are the images in~$\cFF_{2d}$ of the hyperplanes
$H_\beta \subset \cDD_{2d}$ (see e.g.\ \cite[Section~4.4]{mp}). The projection
onto $H^\perp$ given by $v \mapsto \beta - \tfrac1{2d} H$ induces a bijection
of the Noether-Lefschetz divisors $D_{h,a}$ and the Heegner divisor $H_{m,\mu}$
with invariants related by 
\be \label{eq:HgNL}
m \= \frac{a^2}{4d} - (h-1) \, \qquad \mu = a \cdot \omega/2d \bmod L\,.
\ee
(The positivity of~$m$ is guaranteed by the Hodge index theorem.)
We thus use the terms Heegner divisors and Noether-Lefschetz divisors
interchangeably in the K3 case.

\section{Cones of coefficients of modular forms}

Our goal in this section is to show (in Theorem~\ref{thm:modcone} below) 
that the cone generated by the coefficient functionals of Fourier 
expansions of vector-valued modular forms for a lattice~$L$
is rational polyhedral on the space of almost-cusp forms. 
Our main criterion for rational polyhedrality is the following
geometric observation.
\par
\begin{Lemma} \label{le:ratpolycrit}
Suppose that $V$ is a finite-dimensional $\QQ$-vector space
and consider the cone 
\[
\cCC \= \Bigl\{\sum_{n\geq 0} \lambda_n c_n\mid\;\text{all $\lambda_n\in \RR_{\geq 0}$, 
and almost all $\lambda_n$ vanish}\Bigr\}\,\subseteq\, V_\RR = V \otimes \RR
\]
generated by a countable collection of non-zero vectors $(c_n) \subset V$. Suppose
that there exists a codimension one subspace~$S \subset V$ 
and an element $e \in V \setminus S$
with the following properties
\begin{enumerate}
\item[i)] Writing $c_n = \gamma_n e + s$ with $s \in S$, the coefficient 
$\gamma_n$ is positive for all $n$.
\item[ii)] The vectors $c_n$  converge $\RR_{>0}$-projectively to~$e$, 
i.e.\ $c_n/\gamma_n - e\to 0 \in S_\R.$
\item[iii)] Among the $c_n$ there exist elements $c_{n_1},\ldots, c_{n_s}$ such 
that a linear combination $\sum_{i=1}^s \lambda_{i} c_{n_i}$ with all $\lambda_i 
\in \RR_{>0}$ strictly positive lies in~$\langle e \rangle$ and such that
the classes of $c_{n_i} \in V/\langle e \rangle \cong S$ span~$S$.
\end{enumerate}
Then the cone~$\cCC$ is rational polyhedral.
\end{Lemma}
\par
\begin{proof} Due to the first condition the cone lies in the half-space of~$V$
where to $e$-coefficient is positive. It thus suffices to show that the 
convex body~$\cCC_S$ defined as the intersection of~$\cCC$ with the affine 
hyperplane $e + S_\R$ is rational polyhedron. We view $\cCC_S \subset S_\R$ by 
projection along~$e$.
Condition~iii) now implies immediately that this polyhedron $\cCC_S$ contains 
an open neighborhood of zero. Condition~ii) implies that the projections
of~$c_n$ to~$\cCC_S$ converge to zero in~$S_\R$, which is an interior point
of~$\cCC_S$. Consequently, $\cCC_S$ is the convex hull of finitely many
points $\RR_{>0} \cdot c_n \cap (e +S_\R)$. Since these are rational, 
the claim follows.
\end{proof}
\par
We suppose {\em in the remainder of this
section that 
$k \geq 2$,
with $2 k -b^+ +b^- \equiv 4\pmod{8}$,
and that $L = L_1 \oplus U$ splits off a hyperbolic plane~$U$.} 
Here we want to apply Lemma \ref{le:ratpolycrit} to $V = (M_{k,L}^0(\QQ))^\vee$, 
the dual of the space of almost cups forms of half-integral weight~$k$
with rational coefficients. This rationality statement and 
the rationality hypothesis in Lemma~\ref{le:ratpolycrit} is a 
restatement of the fact (\cite{McGraw}) that~$M_{k,L}$ has a basis with
rational coefficients.
\par
The direct sum decomposition~\eqref{eq:ACdecomp}
implies that~$V$ decomposes as 
$$ V \=\langle e \rangle \,\oplus \,S\, $$
where $e$ is defined by $e(E_{k,L}) = -1$ and $e(S_{k,L}) = 0$,
and where $S$ is the subspace of functionals that are zero on~$E_{k,L}$.
We want to apply this lemma to
the vectors~$c_n$ being the coefficient extraction functionals
\[
c_{m,\mu}\,: \, M_{k,L}^0 (\Q)\to \Q, \quad 
f \= \sum_{\mu \in D_L}  \sum_{n\geq 0} a_{m,\mu} q^m e_\mu \mapsto c_{m,\mu}(f)
\= a_{m,\mu}\,.
\]
for $\mu \in D_L$ and $m \in (Q(\mu) + \ZZ) \cap \Q_{> 0}$, i.e.\ the index
set consists of pairs $n=(m,\mu)$.
\par
Condition i) of the lemma is  simply a restatement
of~\eqref{eq:signEIS}. Note that the strict positivity follows from the fact that $L$ splits off a hyperbolic plane over $\Z$. Condition ii) of this lemma is a consequence
of the following proposition.
\par
\begin{Prop} \label{prop:coeffestimates}
Assume $k \geq 5/2$. 
For $\mu\in D_L$ and $m\in Q(\mu)+\Z$ positive, 
 the coefficients $e_{k,L}(m,\mu)$ of the Eisenstein 
series~$E_{k,L}$ 
are negative and satisfy
$$
-e_{k,L}(m,\mu) \,\geq\, C_{k,L} \cdot m^{k-1}
$$
for some positive constant $C_{k,L}$ depending on the weight and the lattice.
\par
If $k=2$ coefficients $e_{k,L}(m,\mu)$ of the Eisenstein series~$E_{k,L}$ 
are negative and for every $\ve>0$ the is a positive constant  $C_{L,\ve}$ depending on 
the lattice such that
$$-e_{2,L}(m,\mu) \,\geq\, C_{L,\ve}\cdot  m^{1-\ve}.
$$
\par
For any $k \geq 2$ the coefficients $a_{m,\mu}$ of any cusp form 
$f = \sum a_{m,\mu} q^m \chi_\mu \in S_{k,L}$
are bounded above in absolute value by 
\be \label{cuspformbound}
| a_{m,\mu}| \, \leq \, C_{f,\ve} \,m^{\frac{k}{2} - \frac{1}{4} + \ve }  
\ee
for some positive constant $C_{f,\ve}$ depending on $f$ and~$\ve$.
\end{Prop}
\par
\begin{proof} 
The last statement is the Weil bound for the coefficients
of cusp forms, see e.g. 
\cite{SarnakModular}, Proposition~1.5.3 and Proposition~1.5.5.
\par
The negativity in the first statement is a consequence of~\eqref{eq:signEIS}
and our congruence condition on~$k$. To prove the lower bound, given the factor
$m^{k-1}$ in both~\eqref{eq:EISeven} and~\eqref{eq:EISodd}, we have to bound
the other terms that depend on~$m$ uniformly from below. 
\par
In the case $b^+$ odd we use the expression in~\eqref{eq:epsod2}.
Lemma~\ref{le:RepLowBd} gives a lower bound for $N_{m,\mu}(p^{w_p} ) / p^{(2k-1)w_p}$
for the finitely many primes dividing~$2N$. For the remaining terms
we use the estimate that for $k \geq 5/2$
\bas
\sum_{g\mid f} \mu(g)\chi_{D'}(g)g^{1/2-k}\sigma_{2-2k}(f/g)
&\,\geq\, 1-\sum_{\substack{g\mid f\\ g>1}} g^{1/2-k}\sigma_{2-2k}(f/g)\\
&\,\geq\, 1-(\zeta(2)-1)\zeta(3) \,\geq\, 1/5.
\eas 
In the case $b^+$ even and $k \geq 3$ we use the expression~\eqref{eq:epseven2}. 
Again, Lemma~\ref{le:RepLowBd} gives a lower bound for the normalized 
representation numbers and together with
\bas
\sigma_{1-k}(d_\mu^2m,\chi_{4D}) \geq 1 - (\zeta(2)-1) > 0
\eas
we obtain a uniform lower bound for $\ve_{m,\mu}$. Finally, in the case $k=2$
we split off the contribution of the divisor~$1$ to $\sigma_{1-k}(d_\mu^2m,\chi_{4D})$
and use the estimate 
$$\sigma_{1-k}(d_\mu^2m,\chi_{4D}) \leq \sigma_1(d_\mu^2m) = O(\log(m))$$
for the remaining terms.
\end{proof}
\par
Because of the direct sum decomposition $V = \langle e \rangle \oplus S$
provided by the conditions in Lemma~\ref{le:ratpolycrit}, its condition~iii) can
be formulated equivalently in terms of the restriction of the
coefficient functionals 
$$\bar{c}_{m,\mu}\, = c_{m,\mu}|_{S_{k,L}(\Q)} : S_{k,L}(\Q) \to \Q$$
to the subspace of cusp forms with rational coefficients. The statement is precisely the content
of the following proposition.
\par

\begin{Prop} 
\label{prop:pos}
There exist 
indices $((m_i,\mu_i))_{i=1}^s$ and real numbers $\lambda_i>0$ such that 
$$ 
\sum_{i=1}^s \lambda_i \,\bar{c}_{m_i,\mu_i} \= 0
$$
in $S_{k,L}(\R)^\vee$ and such that the functionals $(\bar{c}_{m_i,\mu_i})_{i=1}^s$ span~$S_{k,L}(\R)^\vee$.
\end{Prop}
\par
\begin{proof} 
As in \cite{BrPrescribed} we write $L^-$ for the lattice $(L,-Q)$. We identify the Weil represenstation $\rho_{L^-}$ with the dual representation of $\rho_L$. 
The product of a weakly holomorphic modular form~$h$
of weight $2-k$ for the representation $\rho_{L^-}$ 
and any element  $g \in S_{k,L}$ is a weakly holomorphic modular form
of weight~$2$ for  $\Mp_2(\Z)$, i.e.\ a meromorphic differential form on the modular curve $X(1)$. The residue at the cusp $\infty$ of $hg$ vanishes by the residue theorem. 
The idea is to construct a weakly holomorphic modular form~$h$ for 
$\rho_{L^-}$ whose principal part at $\infty$
has non-negative coefficients only 
and whose principal part has sufficiently many non-vanishing terms (that will be the $\lambda_i$ of the proposition) so that the residue
pairing $g \mapsto {\rm Res}(hg)$ involves a spanning system of~$S_{k,L}(\R)^\vee$.
\par
This follows from \cite[Lemma~3.5 and Proposition~3.2]{BrPrescribed}. In fact, 
let 
$$t_\mu = \min \{-Q(\lambda)\,|\; \lambda \in \mu + L, \,\, -Q(\lambda)>0\}\,
\in \tfrac1N \ZZ_{>0} \quad \text{for} \quad \mu \in D_L$$
and let $T = \max\{t_\mu\,|\; \mu \in D_L\}$. 
Choose $B \in \ZZ_{>0}$ sufficiently large such that the weight $k':= 2-k 
+ 12 B>2$ and such that the functionals $\bar{c}_{\mu,\ell}$ for $\ell < B-T$ 
generate~$S_{k,L}(\R)^\vee$.
\par
We claim that the weakly holomorphic modular form
$$ h(\tau) = \Delta(\tau)^{-B} E_{k',L^-}(\tau) \in M^!_{2-k,L^-}$$
has non-negative Fourier coefficients $c_{\mu,\ell}(h)$ and moreover 
\bes
c_{\ell,\mu}(h) >0 \quad \text{for all} \quad  \mu \in D_L, 
\quad \ell \in \ZZ - Q(\mu), \quad \ell \geq T - B\,.
\ees 
To see this, note that the Fourier expansion of $\Delta(\tau)^{-B}$
is $q^{-B}$ times a positive power of the generating function 
$\prod_{j \geq 1}(1-q^j)^{-1}$ of the partition function. 
Consequently, $c_j(\Delta^{-B}) > 0$ for integral 
$j \geq -B$
and zero otherwise. If $\mu \in D_L$ and $\ell \in \ZZ - Q(\mu)$ then
\bes
c_{\ell,\mu}(h) \= c_{\ell - t_\mu} (\Delta^{-B}) \cdot e_{k',L^-}(t_\mu,\mu)
+ \sum_{\ell - t_\mu < j\in \ZZ} c_j(\Delta^{-B}) \cdot e_{k',L^-}(\ell-j,\mu)\,.
\ees
The congruence hypothesis on~$k$ implies that 
$$2k' - b^+(L^-) + b^-(L^-) \= 2k' -b^- + b^+ \,\equiv\, 0 \, (8)$$ 
and hence that by~\eqref{eq:signEIS} all the coefficients $e_{k',L-}(\cdot\,,\,\cdot)$ 
are positive. Consequently, 
the first summand of the right hand side is positive by the definition
of $t_\mu$ and the hypothesis on~$\ell$, and the other summands
are non-negative. This implies the claim.
\par
Let $((m_i,\mu_i))_{i=1}^s$ be some enumeration of the pairs $(m,\mu)$ for $\mu \in D_L$ and $m\in \Z+Q(\mu)$ with
 $0<m < B-T$. By our choice of~$B$, the functionals $\bar{c}_{m_i,\mu_i}$ span~$S_{k,L}(\R)^\vee$. If
we let $\lambda_i = c_{-m_i,\mu_i}(h)$, then the residue theorem
applied to $gh$ implies
$$ \sum_{i=1}^s \lambda_i \,\bar{c}_{m_i,\mu_i}(g) \= 0$$
for any $g \in S_{k, L}$.
\end{proof}
\par
Not only this proof breaks down if the congruence hypothesis on~$k$
is violated, the statement is wrong is this case, as pairing the 
weakly holomorphic form $h$ with the Eisenstein series $E_{k,L}$ bshows, since 
all its coefficients are positive, including the constant term.
\par
We summarize the results of this section to the following statement.
\par
\begin{Thm} \label{thm:modcone}
Let~$L$ be a lattice of signature~$(b^+,b^-)$ that splits off a hyperbolic plane.
We suppose that $k \geq 2$ and $2k - b^+ + b^- \equiv 4 \pmod{8}$.
Then 
the cone~$\cCC$ generated by the coefficient functionals $c_{m,\mu}$ 
on the space of weight~$k$ almost cusp forms $M^0_{k,L}$ for the lattice~$L$
(where $\mu \in D_L$ and $m \in (\ZZ +Q(\mu)) \cap \Q_{> 0}$)   
is a rational polyhedral cone. In particular, the cone~$\cCC$ is
finitely generated.
\end{Thm}
\par

\section{Cones of primitive Heegner divisors} 
\label{sec:IrrNL}

We now translate the results of the previous section into geometric statements.
We suppose for the rest of this paper that  $b^-=2$ and put $k=1+b^+/2$. 
Our motivation for studying cones of coefficient functionals in spaces of modular 
forms comes from the following transport of structure to the rational Picard group 
of orthogonal Shimura varieties. By \cite{BorcherdsGKZ} or 
\cite[Theorem~0.4]{Br-habil}, 
 the map
\be \label{eq:defpsi}
\psi: M_{k,L}^0(\QQ)^\vee \,\to\, \Pic_{\QQ}(\cFF_L(\Gamma)), \qquad 
c_{m,\mu} \mapsto H_{m,\mu}\,,
\ee
sending the coefficient extraction functional $c_{m,\mu}$ to the 
Heegner divisor $H_{m,\mu}$ and the coefficient extraction functional~$-c_{0,0}$.
to the Hodge class $\lambda$, is a homomorphism. We use this to show:
\par
\begin{Thm} \label{thm:redHeeg} 
Let $\Gamma \subseteq \wOL$ be a finite
index subgroup. Suppose that $b^+\geq 3$ and that $L$ splits off a hyperbolic plane.
Then the cone  generated by the (reducible) Heegner divisors $H_{m,\mu}$ 
on $\cFF_L(\Gamma)$ is rational polyhedral. 
\end{Thm}
\par
\begin{proof} 
Since $b^-=2$ and $k=1+b^+/2$, the congruence condition for $k$ in 
Theorem~\ref{thm:modcone} holds. The claim follows from this theorem
and the fact that the image of a rational polyhedral
cone under a linear map is still  rational polyhedral.
\end{proof} 
The map $\psi$ is injective in many situations (e.g.\ if $L$ splits off
two hyperbolic planes, \cite{Br-habil}), we will use this below.
It is moreover surjective (\cite{BLMM}) under 
the hypotheses made here. This implies that the image cone is full-dimensional.
\par
The second goal of this section is to discuss the passage from primitive
to (reducible) Heegner divisors and to prove Theorem~\ref{thm:introEff} 
stated in the introduction, and  the generalization
Theorem~\ref{thm:EffHeeg} below.
\par
\medskip
\paragraph{\bf Primitive Heegner divisors} The Heegner divisors
are in general not irreducible. The divisibility of the defining
lattice element~$v \in L^\vee$ with  $Q(v) = m$ and  $v \equiv \mu 
\mod D_L$ is an obvious invariant distinguishing irreducible components.
Since divisibility is preserved by the action of~$\wOL$, the definition
$$P_{\Delta,\delta} \= \wOL \,\,\backslash \,\,\Biggl( 
\sum_{\substack{L + \delta \ni v \, \text{primitive}, 
\\ Q(v) = \Delta}}  H_v  \Biggr) \,.$$
for $\delta \in D_L$ and $\Delta \in \ZZ + Q(\delta)$ gives well-defined
divisors in $\cFF_L$, called {\em primitive Heegner divisors}. 
By definition (and the fact that any lattice vector
can be written uniquely as a positive multiple of a primitive lattice vector)
\be \label{eq:HinP}
H_{m,\mu} \= \sum_{\substack{r\in\ZZ_{>0}\\ r^2 | m}}\, \sum_{\substack{\delta\in D_L\\ r\delta = \mu}} P_{m/r^2, \delta}\,.
\ee
Here and in the sequel we say that  $r^2\mid m$ with $m \in Q(\mu) + \ZZ$
if there exists $\delta \in D_L$ with $m/r^2 \in Q(\delta) + \ZZ$. We could
drop this condition, since $P_{m/r^2, \delta}$ is empty otherwise.
\par
The converse to~\eqref{eq:HinP} follows from a variant of M\"obius inversion.
\par
\begin{Lemma} \label{le:IrrHeeg}
The primitive divisors $P_{\Delta, \delta}$ can be written in terms of the
Heegner divisors $H_{m,\sigma}$ as 
\be \label{eq:IrrHeeg}
P_{\Delta,\delta} \= \sum_{\substack{r\in \ZZ_{>0}\\ r^2 | \Delta}} \mu(r) \sum_{\substack{\sigma\in D_L\\ r\sigma  = \delta}} 
H_{\Delta/r^2, \sigma}\,.
\ee 
\end{Lemma}
\par
\begin{proof} The statement for $\Delta$ without quadratic divisors
is obvious and from the definition in~\eqref{eq:HinP} and 
induction on the number of quadratic divisors we obtain
\ba P_{\Delta,\delta} &\=
 H_{\Delta, \delta} - \sum_{1 \neq s :  s^2 | \Delta} \,
\sum_{\tau s = \delta}
\Biggl(\,
\sum_{t: t^2|\Delta/s^2} \mu(t) \sum_{\sigma t = \tau}  
H_{\Delta/ s^2t^2, \sigma} \Biggr) \\
&\=
H_{\Delta, \delta} - \sum_{1 \neq s :  s^2 | \Delta} 
\Biggl(\,
\sum_{t: t^2|\Delta/s^2} \mu(t) \sum_{\tau s = \delta} \sum_ {\sigma t = \tau} 
H_{\Delta/ s^2t^2, \sigma} \Biggr)\,.\\  
\ea
We can group the interior double sum as a single sum over all~$\sigma$
with $\sigma \cdot st = \delta$. We let~$r = st$ and consider the
summands contributing to $H_{\Delta/r^2,\sigma}$. It remains to show
that for given $r$ the exterior double sum including the factor~$\mu(t)$
adds up to $\mu(r)$. If some prime divides $r$ more then once, the
claimed contribution follows since $\sum_{I \subseteq P} (-1)^{|I|} = 0 = \mu(h)$ 
for any finite set (of primes)~$P$. In the remaining cases, one summand
is missing in the subset summation since $1 \neq s$, and with the global
minus sign we obtain the coefficient $\mu(h)$ we want.
\end{proof}
\par
\medskip
\paragraph{\bf Primitive Noether-Lefschetz divisors}
In the case of $K3$-lattice $L = L_{2d}$ the decomposition of Heegner 
divisors into irreducible components can also be described by
a geometric decomposition of Noether-Lefschetz divisors.
The Picard group of a generic member of the Noether-Lefschetz divisors
contains a rank two lattice $\Lambda$ with signature $(1,1)$. Conversely, the 
intersection matrix $\Lambda$ with a distinguished element $H$ with $Q(H) = 2d$
has the intersection form 
$$M_{\Lambda} = \mat {2d} y y {2x} $$
with respect to some basis $\{H,\beta\}$ starting with~$H$.
The discriminant $\Delta(\Lambda) = \det(M_\Lambda) \in \ZZ$ 
and the coset $\delta = y \mod 2d$ are invariants
of such a lattice and it is easy to show that the pair $(\Delta, \delta)$
is a classifying invariant of such rank two lattices.
We now define the primitive Noether-Lefschetz divisors 
$\cPP_{\Delta, \delta}$ to be the closure of the locus of K3-surfaces that have a 
sublattice $\Lambda \subset L_{K3}$ of signature $(1,1)$, containing~$H$
and we provide them with multiplicity one or two depending on $2\delta \neq 0$
or not modulo~$2d$. The similarity in notation to primitive Heegner divisors 
is intentional, since we claim that
\be
\cPP_{\Delta, \delta}  \= P_{\Delta/4d, \delta (\omega/2d)}\,.
\ee
This can be seen from the definitions, tracing the definitions 
along the bijection $v \mapsto \beta - \tfrac{1}{2d}H$ between
Heegner and Noether-Lefschetz divisors given along with
equation~\eqref{eq:HgNL}.
\par
\medskip
\paragraph{\bf The main result on ${\rm Eff}^{{\rm H}}(\cFF_{L})$}
We suppose 
in the remainder of this section that $b^+\geq 3$ and that 
$L = L_0 \oplus U^{\oplus 2}$ splits off {\bf two} copies of 
a hyperbolic plane~$U$. We put $k=1+b^+/2$ and $\Gamma = \wOL$, i.e.\ we work on 
the Shimura varieties $\cFF_L$.
\par
\begin{Lemma}
Under these conditions the primitive Heegner divisors $P_{\Delta, \delta}$ are 
irreducible if $2\delta \neq 0 \in D_L$. If $2\delta = 0 \in D_L$
then $P_{\Delta, \delta}$ is (if non-empty) an irreducible divisor
with multiplicity two.
\end{Lemma}
\par
\begin{proof} The multiplicity two stems from the fact that~$v$
and $-v$ give the same divisor if $2\delta = 0 \in D_L$. It remains
to show that any two primitive elements in~$L^\vee$ with the same
norm and the same $D_L$-coset lie in the same $\wOL$-orbit. This
can be done using Eichler-transformations, see Lemma~4.4 in~\cite{FrHe}.
\end{proof}
\par
The primitive Heegner divisors in our main result are thus irreducible. 
\par
\begin{Thm} \label{thm:EffHeeg}
The cone  ${\rm Eff}^{{\rm H}}(\cFF_{L})$ generated by the primitive  
Heegner divisors $P_{\Delta, \delta}$ is rational polyhedral. In particular it is 
finitely generated.
\end{Thm}
\par
\smallskip
We prove this theorem in the same way as Theorem~\ref{thm:modcone}, using a
refinement of the estimate in Proposition~\ref{prop:coeffestimates}. 
Together with the main observation of the proof of Lemma~\ref{le:ratpolycrit}, the
following proposition implies that the vectors corresponding to primitive 
Heegner divisors still converge to an interior point of the cone
$\cCC = {\rm Eff}^{{\rm H}}(\cFF_{L})$. Theorem~\ref{thm:introEff} 
and Theorem~\ref{thm:EffHeeg} follow immediately from the following statement.
Let $\varphi = \psi^{-1}$ be the inverse of the map~$\psi$ defined 
in~\eqref{eq:defpsi}.
\par
\begin{Prop} Under the identification $\varphi$ any infinite sequence
of pairwise different primitive Heegner divisor $P_{\Delta, \delta}$
converges  $\RR_{>0}$-projectively to the functional~$e$.
\end{Prop} 
\par
\begin{proof} As in Proposition~\ref{prop:coeffestimates}, we show that there 
is a constant $C>0$ and for any $\ve>0$ and any cusp form~$f$ of weight~$k$ 
there is a constant $C_{f,\ve}>0$ 
such that the bounds 
$$\varphi(P_{\Delta, \delta})(E_{k,L}) \,\geq \, C \cdot \Delta^{k-1} $$
for the coefficients functional evaluated at Eisenstein series, and 
$$
|\varphi(P_{\Delta, \delta})(f)| \,\leq\, C_{f,\ve} \cdot \Delta^{\frac{k}2 - 
\frac{1}{4} + \ve} 
\quad \text{if} \quad k \in \ZZ$$
hold for any $\delta \in D_L$ and any positive $\Delta \in Q(\delta) + \ZZ$.
\par
The claim about cusp forms follows immediately from~\eqref{cuspformbound}, 
since the number of summands in~\eqref{eq:IrrHeeg} contributing
to $P_{\Delta,\delta}$ is at most $D_L$ times the number of square divisors
of $\Delta$, which is $O(\Delta^\ve)$ for any~$\ve$.
\par
In order to estimate the Eisenstein series contribution,we define
$K_r \subseteq D_L$ to be the kernel of the multiplication by~$r$ and
we observe that $1 \leq |K_r| \leq r^{2k-4}$, since the lattice~$L_0$ 
is of rank~$2k-4$. 
\par
In the case $b^+$~odd, we deduce from~\eqref{eq:IrrHeeg} and the 
formula~\eqref{eq:EISodd}
for the coefficients that $\varphi(P_{\Delta,\delta})$ evaluated at 
the Eisenstein series
is equal to $\Delta^{k-1}$ times some constants independent of~$\Delta$ times
\be \label{eq:ef}
Q(\Delta,\delta) \= \sum_{r: r^2|\Delta} \mu(r) \sum_{\mu r =\delta} 
\frac{1}{r^{2(k-1)}}  \ve_{\Delta/r^2,\mu}\,,
\ee
where $\ve_{m,\mu}$ was defined in~\eqref{eq:epsod1}.
We want to show the that there is some~$C>0$ such that $Q(\Delta,\delta) > C$ 
for all~$\Delta$. By definition of the M\"obius function, 
$Q(\Delta,\delta)$ is greater or equal to $\ve_{\Delta,\delta}$ (stemming 
from $r=1$) minus the sum over subsets~$P$ of odd cardinality 
of the set of prime divisors of~$\Delta$. In order to estimate these
negative contributions from above, we compare $\ve_{\Delta,\delta}$
with $\ve_{\Delta/r^2,\mu}$ using~\eqref{eq:epsod1}. First we remark that 
we can arbitrarily enlarge (by Theorem~7 in \cite{BrKuss}) the set of primes over which the product runs.
Hence we can suppose that the product runs over the same set of primes when
computing $\ve_{\Delta,\delta}$ and $\ve_{\Delta/r^2,\mu}$. 
The terms $1-p^{1-2k}$
obviously cancel and we claim that the same holds for the terms
$1-\chi_{D'}(p)p^{1/2-k}$. Here $D'=2(-1)^{(b^++b^-+1)/2}\Delta|D_L|$ and
the corresponding discriminant associated with $\Delta/r^2$
is by definition $\widetilde{D}' = D'/r_0^2$, where we have written $r =r_0 r_1$
with $r_1$ the largest factor in~$r$ coprime to~$2N$. Said differently, 
$\widetilde{D}'$ and~$D'$ differ only in prime factors~$p$ dividing~$2N$
and for those $\chi_{D'}(p) = 0 = \chi_{\widetilde{D}'}(p)$ since $2|D_L|$ 
divides~$\widetilde{D}'$.
\par
As quotient of the remaining factors we obtain 
\ba \label{eq:vequot}
\frac{\ve_{\Delta/r^2,\mu}}{\ve_{\Delta,\delta}}&= \prod_{p\mid r} \frac{p^{w_p
(\Delta/r^2,\mu)(1-2k)}N_{\Delta/r^2,\mu}(p^{w_p
(\Delta/r^2,\mu)})}{p^{w_p(\Delta,\delta)(1-2k)}N_{
\Delta,\delta}(p^{w_p
(\Delta,\delta)})}\\
&= \prod_{p\mid r} \frac{p^{w_p
(\Delta,\delta)(1-2k)}N_{\Delta/r^2,\mu}(p^{w_p
(\Delta,\delta)})}{p^{w_p(\Delta,\delta)(1-2k)}N_{
\Delta,\delta}(p^{w_p
(\Delta,\delta)})}
.
\ea
According to Lemma~\ref{le:RepLowBd2}, we get the bound
\be \label{eq:veest}
\frac{\ve_{\Delta/r^2,\mu}}{\ve_{\Delta,\delta}}\leq \prod_{p\mid r} 
\frac{1+p^{-1}}{1-p^{-2}}= \prod_{p\mid r} 
\frac{1}{1-p^{-1}}.
\ee
Using $|K_r| \leq r^{2k-4}$, we obtain
\ba
\frac{Q(\Delta,\delta)}{\ve_{\Delta,\delta}} &\,\geq\, 1
\,-\, \sum_{|P| \,{\rm odd}} \prod_{p \in P} 
\frac{1}{p(p-1)} \\
& \,\geq\, 
1-\frac12 \Bigl(\prod_{p\,{\rm prime}} \Bigl(1+\frac{1}{p(p-1)}\Bigr)  
- \prod_{p\,{\rm prime}} \Bigl(1-\frac{1}{p(p-1)}\Bigr) \Bigr) .
\ea
The Euler products appearing on the right hand side are known as 
Landau's totient constant
\[
\prod_{p\,{\rm prime}} \Bigl(1+\frac{1}{p(p-1)}\Bigr)= 
\frac{315}{2\pi^4}\zeta(3) = 1.943596...
\]
and Artin's constant
\[
\prod_{p\,{\rm prime}} \Bigl(1-\frac{1}{p(p-1)}\Bigr)= 0.373955...,
\]
respectively. Inserting the numerical values, we see that 
\ba
\frac{Q(\Delta,\delta)}{\ve_{\Delta,\delta}} &\,\geq\, 
0.215179....\,>\,0,
\ea
which can be rigorously proven to be positive by standard
remainder term estimaes for zeta-functions.
\par
The case $b^+$ even is similar, but easier. Again we need to estimate
$\frac{\ve_{\Delta/r^2,\mu}}{\ve_{\Delta,\delta}}$ from~\eqref{eq:epseven1} uniformly
from below. By Theorem~7 in \cite{BrKuss} we may again suppose that
the product runs over the same set of primes for~$\Delta$ and~$\Delta/r^2$.
This time, the discimrinant involved in $1-\chi_D(p)p^{-k}$
does not depend on~$m$. Hence the corresponding factors cancel. The
remaining expression is the same as in~\eqref{eq:vequot}
and can be estimated as in~\eqref{eq:veest} above. 
\end{proof}

\section{The case of Hilbert modular surfaces} \label{sec:HMS}

Let $F=\Q(\sqrt{D})$ be a real quadratic field of discriminant $D>0$ with ring of integers $\calO_F$. Denote the conjugation in $F$ by $\nu\mapsto \nu'$. The Hilbert modular surface associated with $F$  and an ideal $\frakb\subset \cOO_F$ is the quotient 
$X_{F,\frakb} = \SL(\cOO_F \oplus \frakb)\backslash\HH^2$,
see e.g.\ \cite{vandergeer88} for a textbook reference. 
In this section we show that the results of the preceding sections 
partially also apply to Hilbert modular surfaces.
%
\par
We first remark that the cone generated by Heegner-divisors (also
called Hirze\-bruch-Zagier cycles in this case) is no longer full dimensional 
on $X_{F,\frakb}$, contrary to the case of $b^+ \geq 3$ (see \cite{BLMM}).
First, the two foliations on Hilbert modular surfaces
define two line bundles $\cLL_1$ and $\cLL_2$ on $X_{F,\frakb}$.
Heegner-Divisors always lie in the subspace whose intersection with
$\cLL_1 \otimes \cLL_2^{-1}$ is zero. Second, even in this subspace, the cone
is not always full-dimensional, see \cite{HLR}.
\par
Hilbert modular surfaces nevertheless fall into the scope of the
preceding sections. In fact, let $B=\norm(\frakb)$, and consider the lattice 
\[
L_\frakb= \left\{ \zxz{x}{\nu'}{\nu}{y}:\; x\in \Z,\, y\in B\Z, \, \nu\in \frakb\right\}
\]
with the integral quadratic form $Q(X)= -\frac{1}{B}\det(X)$. The dual lattice of $L_\frakb$ is given by 
\[
L_\frakb^\vee= \left\{ \zxz{x}{\nu'}{\nu}{y}:\; x\in \Z,\, y\in B\Z, \, \nu\in \frakd_F^{-1}\frakb\right\},
\]
where $\frakd_F\subset\calO_F$ is the different ideal. In particular, we have $L_\frakb^\vee/L_\frakb\cong \frakd_F^{-1}/\calO_F$.
The Hilbert modular group $\SL(\cOO_F \oplus \frakb)$ acts on $L_\frakb$ by 
\[
(g,X)\mapsto g X\,{}^tg'
\]
for $g\in \SL(\cOO_F \oplus \frakb)$ and $X\in L_\frakb$. This action preserves the quadratic form $Q$, and according to \cite[Section~2.7]{Br-123} it induces an isomorphism 
\[
\SL(\cOO_F \oplus \frakb)\cong\operatorname{Spin}(L_\frakb) .
\]
On the other hand, according to \cite[Lemma~2.6]{Madapusi}, the image of the spin group of a lattice in the orthogonal group is given by the intersection of the stable special orthogonal group with the subgroup of elements with positive spinor norm, that is, 
\[
\operatorname{Spin}(L_\frakb)/\{\pm 1\} \cong \widetilde\SO^+(L_\frakb).
\]
Consequently, we have 
$\SL(\cOO_F \oplus \frakb)/\{\pm 1\} \cong \widetilde\SO^+(L_\frakb)$, and 
\[
X_{F,\frakb}\cong \calF_{L_\frakb}(\widetilde\SO^+(L_\frakb)).
\]
The explicit identification of $X_{F,\frakb}$ with the orthogonal Shimura variety on the right hand side is given by \cite[Equation (2.33)]{Br-123}.
The Heegner divisors on the right hand side can be identified with Hirzebruch-Zagier divisors on $X_{F,\frakb}$. 

To describe the symmetric Hilbert modular group, we consider the vector $\lambda=\kzxz{1}{0}{0}{-B}\in L_\frakb$ with $Q(\lambda)=1$. The reflection $\tau_\lambda\in \Orth(L_\frakb)$ taking $\lambda$ to its negative and fixing its orthogonal complement belongs to $\widetilde \Orth^+(L_\frakb)$ and has determinant $-1$.
On $\H\times \H$ it induces the transformation 
\begin{align}
\label{eq:syminv}
(z_1,z_2)\mapsto \left(-\frac{1}{Bz_2},-\frac{1}{Bz_1}\right).
\end{align}
Hence, the projective symmetric Hilbert modular group is isomorphic to $\widetilde \Orth^+(L_\frakb)$, and  the corresponding symmetric Hilbert modular suface is given by  
\[
X_{F,\frakb}^{\rm symm}\cong \calF_{L_\frakb}(\widetilde\Orth^+(L_\frakb)).
\]
Since 
$L_\frakb$ splits one hyperbolic plane over $\Z$, we may apply Theorem \ref{thm:redHeeg} in this situation.
\par
\begin{Cor}
The cone generated by the (reducible) Hirzebruch-Zagier cycles on the 
Hilbert modular suface $X_{F,\frakb}$ is rational polyhedral. The same
statement holds on the symmetric Hilbert modular surface $X_{F,\frakb}^{\rm symm}$.
\end{Cor}
\par
It seems quite plausible that the rational polyhedrality can be extended
in the case of Hilbert modular sufaces to the cone generated by the
irreducible components of Hirzebruch-Zagier cycles. The description
is more complicated than in the case when $L$ splits off two hyperbolic
planes. It has been given in many cases by Hirzebruch, Franke and
Hausmann, see e.g.\ 
the survey and references in \cite[Section~5.2]{MoeZag}
or \cite[Section~5.3]{vandergeer88}.

\bibliography{my}
\end{document}